\documentclass[10pt]{amsart}

\usepackage{hyperref,amsmath,amsthm,amsfonts,amscd,flafter,epsf,amssymb,diagrams}
\usepackage{epsfig}
\usepackage{amsfonts}
\usepackage{amssymb}
\usepackage{amsmath}
\usepackage{amsthm}
\usepackage{latexsym}
\usepackage{amscd}
\usepackage{flafter}
\usepackage{epsf}
\input{diagrams.sty}
\newcommand{\ra}{\rightarrow}
\newcommand{\xra}{\xrightarrow}
\newcommand{\Sig}{|Sigma}
\newcommand{\thetabar}{{\overline{\theta}}}
\newcommand{\phibar}{{\overline{\phi}}}
\newcommand{\psibar}{{\overline{\psi}}}
\newcommand{\Psibar}{{\overline{\Psi}}}
\newcommand{\Phibar}{{\overline{\Phi}}}
\newcommand{\rrr}{\mathfrak{r}}
\newcommand{\sss}{\mathfrak{s}}
\newcommand{\xxx}{\mathfrak{x}}
\newcommand{\yyy}{\mathfrak{y}}
\newcommand{\zzz}{\mathfrak{z}}
\newcommand{\ppp}{\mathfrak{p}}

\newcommand{\xla}{\xleftarrow}

\newcommand\Dual{\mathcal D}
\newcommand\Duality\Dual

\newcommand\relspinc{\underline{\spinc}}

\newcommand\ModSphere{\ModFlow\left({\mathbb S}\longrightarrow
\Sym^{g-1}(\Sigma_{1})\times \Sym^2(\Sigma_{2})\right)}
\newcommand\ModSpheres\ModSphere

\newcommand\UnparModSp{\widehat \ModSp}
\newcommand\UnparModFlow\UnparModSp

\newcommand{\spinc}{\mathfrak s}

\newcommand\ModMaps{\mathcal M}
\newcommand\ModSp\ModMaps

\newcommand\alphas{\mbox{\boldmath$\alpha$}}

\newcommand\betas{\mbox{\boldmath$\beta$}}

\newcommand\spincrel\relspinc

\newtheorem{thm}{Theorem}[section]

\def\endproof{\relax\ifmmode\expandafter\endproofmath\else
  \unskip\nobreak\hfil\penalty50\hskip.75em\hbox{}\nobreak\hfil\bull
  {\parfillskip=0pt \finalhyphendemerits=0 \bigbreak}\fi}
\def\endproofmath$${\eqno\bull$$\bigbreak}
\def\bull{\vbox{\hrule\hbox{\vrule\kern3pt\vbox{\kern6pt}\kern3pt\vrule}\hrule}}

\newcommand{\Z}{\mathbb{Z}}

\newcommand{\ModSWfour}{\mathcal{M}}
\newcommand{\ModFlow}{\ModSWfour}

\newcommand\abuts\Rightarrow
\newcommand\Sym{\mathrm{Sym}}

\newcommand{\Hbb}{{\mathbb{H}}}

\begin{document}

\title[A combinatorial approach to surgery]{A combinatorial approach to surgery formulas in Heegaard Floer homology}%
\author{Eaman Eftekhary}%
\address{School of Mathematics, Institute for Studies in Theoretical Physics and Mathematics (IPM),
P. O. Box 19395-5746, Tehran, Iran}%
\email{eaman@ipm.ir}

\thanks{The authors is partially supported by a NSF grant}%
\keywords{Heegaard Floer homology, surgery formula}%

\maketitle
\begin{abstract}
Using the combinatorial approach to Heegaard Floer homology, we obtain
a relatively easy formula for  computation of $\widehat{\text{HF}}(Y_{p/q}(K),\Z_2)$,
where $Y_{p/q}(K)$ is the three-manifold obtained by $p/q$-surgery on a
knot $K$ inside a homology sphere $Y$.
\end{abstract}
\section{Introduction}

In \cite{Ef-splicing}, the author used the combinatorial description
of Heegaard Floer homology (see \cite{SW,MOS} and also \cite{MOST})
to obtain a gluing formula for Heegaard Floer homology, when two bordered three-manifolds
with torus boundary are glued together.
In this paper, we apply the gluing formula of \cite{Ef-splicing} to the
especial case of rational surgeries on the knots inside homology spheres,
and after simplifications, we derive a relatively easy formula for the
Heegaard Floer homology of these three-manifolds. Let $K$ be a knot inside
the homology sphere $Y$. We may remove a tubular neighborhood of $K$ and
glue it back to obtain the three-manifold $Y_{p/q}=Y_{p/q}(K)$, which is the
result of $p/q$-surgery on $K$. The core of the solid torus, which is the
tubular neighborhood of $K$, will represent a knot in $Y_{p/q}$ which will be
denoted by $K_{p/q}$. We may denote $(Y,K)$ by $(Y_\infty,K_\infty)$, as
an extension of the above notation. Let $\Hbb_\bullet(K)$ be the Heegaard
Floer homology group $\widehat{\text{HFK}}(Y_\bullet,K_\bullet)$ for
$\bullet\in \mathbb{Q}\cup\{\infty\}$. Note that $\widehat{\text{HFK}}$ is
defined for knots inside rational homology spheres (see \cite{OS-Qsurgery}), and that
$\Hbb_0(K)=\widehat{\text{HFL}}(Y,K)$ is the longitude Floer homology
of $K$ from \cite{Ef-Whitehead}. In all these cases, we choose the
coefficient ring to be $\Z/2\Z$.\\

\begin{figure}
\mbox{\vbox{\epsfbox{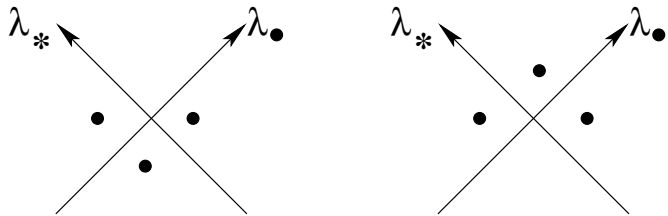}}} \caption{\label{fig:punctures}
{For defining chain maps between $C_\bullet(K)$ and $C_\star(K)$, the
punctures around the intersection point of $\lambda_\bullet$ and $\lambda_\star$
should be chosen as illustrated in the above diagrams.
}}
\end{figure}

If we choose a Heegaard diagram for $Y-K$ and let $\lambda_\bullet$ denote
a longitude which has framing $\bullet\in \Z\cup\{\infty\}$ (with $\lambda_\infty=\mu$
the meridian for $K$), one may observe that the pairs $(\lambda_\infty,\lambda_1)$ and
$(\lambda_1,\lambda_0)$ have a single intersection point in the Heegaard diagram.
Let $(\bullet,\star)\in\{(\infty,1),(1,0)\}$ correspond to either of these pairs.
There are four quadrants around the intersection point of $\lambda_\star$ and
$\lambda_\bullet$. If we puncture three of these quadrants and consider the
corresponding holomorphic triangle map, we obtain an induced map $\Hbb_\bullet\ra
\Hbb_\star$. If the punctures are chosen as in figure~\ref{fig:punctures}, the
result would be  two maps
$\phi,\phibar:\Hbb_\infty(K)\ra \Hbb_1(K)$ and two other maps
$\psi,\psibar:\Hbb_1(K)\ra \Hbb_0(K)$ so that the following
two sequences are exact:
\begin{displaymath}
\begin{split}
&\Hbb_\infty(K) \xra{\phi}\Hbb_1(K)\xra{\psibar}\Hbb_0(K),\ \ \ \text{and}\\
&\Hbb_\infty(K) \xra{\phibar}\Hbb_1(K)\xra{\psi}\Hbb_0(K).\\
\end{split}
\end{displaymath}
The homology of the mapping cones of $\phi$ (or $\phibar$) and $\psi$ (or $\psibar$)
are $\Hbb_0(K)$ and $\Hbb_\infty(K)$ respectively (see \cite{Ef-splicing}).
With the above notation fixed, we prove the following surgery formula:
\begin{thm}\label{thm:main}
Let $K$ be a knot in a homology sphere $Y$ and let the complexes $\Hbb_\bullet=\Hbb_\bullet(K)$, $\bullet
\in \{\infty,1,0\}$ and the maps $\phi,\phibar,\psi,\psibar$ between them be as above.
The homology of $Y_{p/q}(K)$, the manifold obtained by $\frac{p}{q}$-surgery on $K$
(for positive integers $p,q$ with $(p,q)=1$),
may be obtained as the homology of the complex $(\Hbb,d)$ where
$$\Hbb=\Big(\bigoplus_{i=1}^q\Hbb_\infty(i)\Big)\oplus\Big(\bigoplus_{i=1}^{p+q}\Hbb_1(i)\Big)
\oplus\Big(\bigoplus_{i=1}^p\Hbb_0(i)\Big),$$
where each $\Hbb_\bullet(i)$ is a copy of $\Hbb_\bullet$. Moreover, the differential $d$ is the sum of
the following maps
\begin{displaymath}
 \begin{split}
  \phi^i:\Hbb_\infty(i)\ra \Hbb_1(i),\ \ \ \ \ \ \ \ &\phibar^i:\Hbb_\infty(i)\ra \Hbb_1(i+p),\ \ \ i=1,2,...,q\\
  \psi^j:\Hbb_1(j+q)\ra \Hbb_0(j), \ \ \ &\phibar^j:\Hbb_1(j)\ra \Hbb_0(j),\ \ \ \ \ \ \ \ j=1,2,...,p,\\
 \end{split}
\end{displaymath}
where $\phi^i$ is the map $\phi$ corresponding to the copy $\Hbb_\infty(i)$ of $\Hbb_\infty$, etc..
\end{thm}

Surgery formulas for Heegaard Floer homology were first studied by
Ozsv\'ath and Szab\'o in \cite{OS-surgery} and \cite{OS-Qsurgery}.
Their work was generalized
to a complete description of the quasi-isomorphism type of the
Heegaard Floer complex associated with
$(Y_n(K),K_n)$ in terms of $\text{CFK}^\infty(Y,K)$ by the author in
\cite{Ef-surgery}. The gluing formulas here are, however, simpler for
actual computations. They are also used in the proof of the main theorem of
\cite{Ef-incompressible}, which states that a prime homology sphere with trivial
Heegaard Floer homology can not contain an incompressible torus.\\

{\bf{Acknowledgement.}} I would like to thank Matt Hedden for useful discussions which 
resulted in correcting a mistake in an earlier version of this paper. 

\section{Heegaard diagrams for surgery on a knot}
The Heegaard diagram for a knot $K$ in a three-manifold $Y$
is a surface $\Sigma$ of genus $g$ together with a set
$\alphas=\{\alpha_1,...,\alpha_g\}$ of $g$ homologically linearly independent
pairwise disjoint simple closed curves on $\Sigma$ and another set $\betas_0=
\{\beta_1,...,\beta_{g-1}\}$ of such curves together with two especial curves
$\mu$ and $\lambda$. The curves $\mu$ and $\lambda$ intersect each-other once,
and transversely, in a single point $z$,
while both of them are disjoint from all the curves in $\betas_0$. They are
chosen so that $(\Sig,\alphas,\betas_0)$ is a Heegaard diagram for
$Y-\text{nd}(K)$, and such that $\mu$ represents the meridian
of the knot $K$. The curve $\lambda$ corresponds to the unique simple closed curve
in the boundary of $\text{nd}(K)$ which is trivial in the first homology of
$Y$.
We may assume that  the complement of these curves in $\Sigma$ will be a union of
regions of a very especial type; one of the regions may be
arbitrarily complicated (we call it the \emph{bad} region).
The rest of them consist of three triangles, a number
of rectangles, and a number of bi-gons. Each triangle will have $z$ as a
corner. The bad region will also have $z$ as a corner.\\

Associated with this setup,  two chain complexes $M(K)$ and $L(K)$
are constructed. The complex $M(K)$ is the Heegaard chain complex
associated with the Heegaard diagram
$$(\Sig,\alphas,\betas_0\cup\{\mu\},u,v),$$
where $u$ and $v$ are two marked points on the two sides of $\mu$. It is the
familiar complex $\widehat{\text{CFK}}(Y,K)$. Similarly, we may use $\lambda$
instead of $\mu$ to obtain the complex $L(K)$, which is the chain complex
$\widehat{\text{CFL}}(Y,K)$ from \cite{Ef-Whitehead}. In \cite{Ef-splicing}
we constructed four
chain maps in this situation $\Psi_{i}:L(K)\ra M(K)$, $i=1,2,3$ and
$\Phi:M(K)\ra L(K)$. They satisfied $\Psi_2\circ \Phi \circ \Psi_1=\Psi_3$.
The maps $\Psi_1,\Phi$ and $\Psi_2$ correspond to counting triangles which use
the first, the second and the third quadrant around the intersection point $z$ respectively
(in the clockwise order) which are not part of the bad region.
\\

Suppose that the knots $(Y_i,K_i),\ i=1,2$ are given, and that $Y$ is the three-manifold
obtained by splicing the complement of $K_1$ and the complement of $K_2$.
Let $M(K_i)$ and $L(K_i)$ and $\Psi^i_j, j=1,2,3$ and $\Phi^i$ be the
complexes and the maps associated with $(Y_i,K_i)$ for $i=1,2$. Let
$M(Y)=M(K_1)\otimes M(K_2)$ and $L(Y)=L(K_1)\otimes L(K_2)$ and
construct $C(Y)$ as  the complex $L(Y)\oplus M(Y)$.
Let $d_C$ be the differential of this complex. Define the maps
$\Phibar=\Phi^1\otimes \Phi^2:M(Y)\ra L(Y)$ and
$$\Psibar=\Psi^1_1\otimes \Psi^2_2+\Psi^1_2\otimes \Psi^2_1+\Psi^1_3\otimes
\Psi^2_3:L(Y)\ra M(Y),$$ and set $D_C=d_C+\Phibar+\Psibar$.
The following was shown in \cite{Ef-splicing}:

\begin{thm}
In the above situation and with the above notation, $(C(Y),D_C)$ is a
chain complex. The homology of this chain complex gives the Heegaard
Floer homology group $\widehat{\text{\emph{HF}}}(Y)$.
\end{thm}

\section{A diagram in the lens space}

\begin{figure}
\mbox{\vbox{\epsfbox{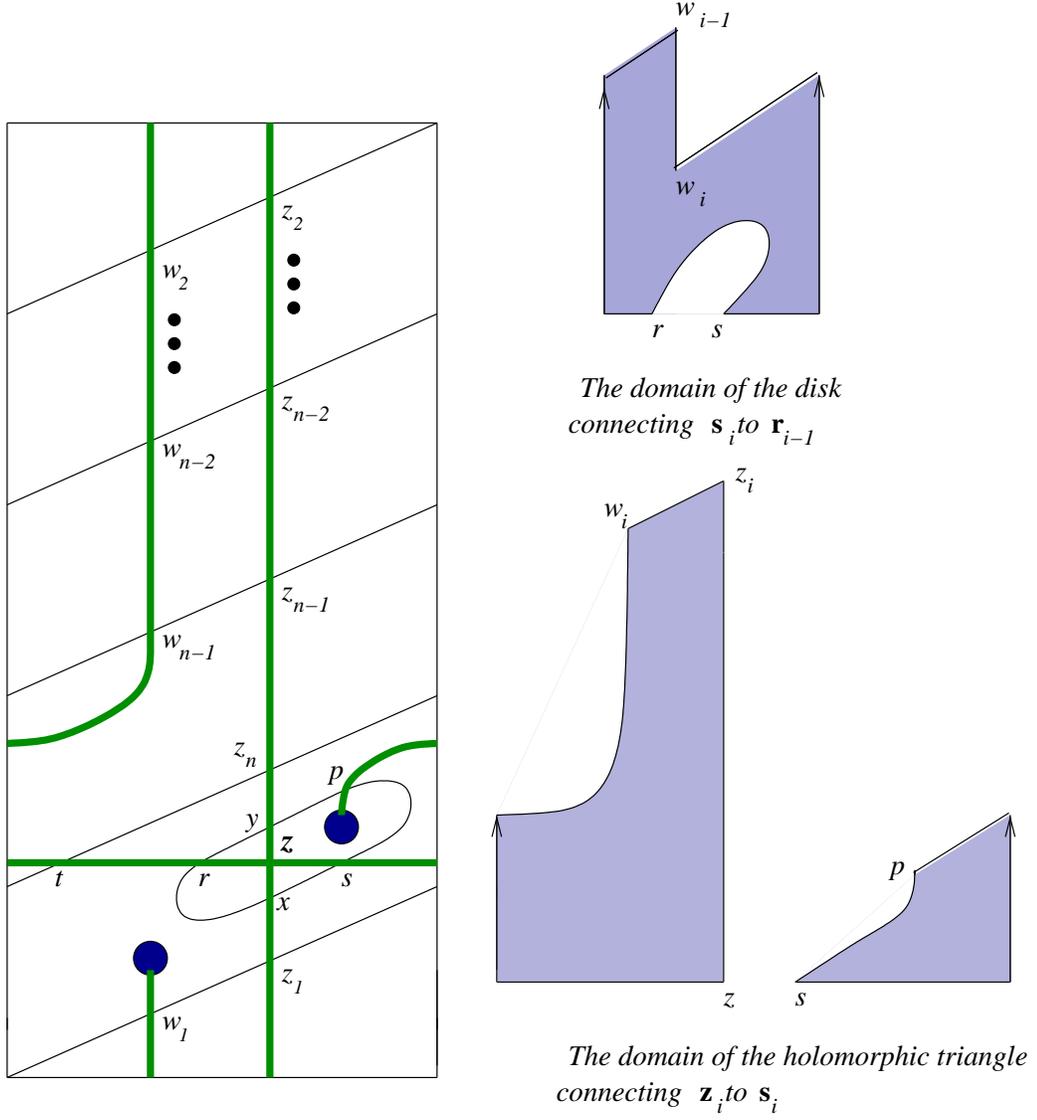}}} \caption{\label{fig:HD}
{The especial Heegaard diagram which should be used for $n$-surgery
on knots inside homology spheres is illustrated.
}}
\end{figure}
 Let $n$ be a given
positive integer and consider the Heegaard diagram $H_n$
illustrated in figure~\ref{fig:HD}.
Let $M(n)$ and $L(n)$ be the chain complexes associated with this diagram.
In this section, we will discuss these chain complexes, together with the
corresponding maps $\Psi_i(n):L(n)\ra M(n)$ and $\Phi(n):M(n)\ra L(n)$.

Each generator of the complex $L(n)$ consists of two intersection points in the
Heegaard diagram $H_n$. One of the intersection points $r,s$ or $t$ on the
curve $\lambda_n$ has to be chosen. If $r$ or $s$ is chosen, we may complete
this choice to a generator for $L(n)$ by choosing one of $w_i,\ i=1,...,n-1$.
Denote the pair $\{r,w_i\}$ by $\rrr_i$. Similarly, if we use $s$ instead of $r$
we obtain $\sss_i,\ i=1,...,n-1$. The intersection point $t$ may be completed
 to a generator only in one way, namely by pairing it with $p$. Denote this generator
by $\ppp$. It is not hard to see in the picture that there is a disk, with a domain that
is topologically a cylinder, which goes from $\sss_{i+1}$ to $\rrr_i$ for $i=1,...,n-1$.  Furthermore, there is a differential going from
$\ppp$ to $\rrr_{n-1}$. As a result, the differential of $L(n)$ may be
described as follows:
\begin{displaymath}
\begin{split}
&\ell_n(\sss_{i+1})=\rrr_i,\ \ \ \ \ i=1,2,...,n-1\\
&\ell_n(\ppp)=\rrr_{n-1}\\
&\ell_n(\xxx)=0, \ \ \ \ \ \ \ \ \text{if $\xxx$ is any other generator.}
\end{split}
\end{displaymath}
Similarly, a generator of the complex $M(n)$ will contain one of the three
intersection points $x,y$ or $p$. The intersection points $x$ and $y$ may
be completed to a generator by adding one of the intersection points $w_i,\
i=1,...,n-1$. This way we get the generators $\xxx_i$ and $\yyy_i$ for
$i=1,...,n-1$. The intersection point $p$ may be completed to a generator
by adding one of the intersection points $z_i,\ i=1,...,n$. The corresponding
generators will be denoted by $\zzz_i,\ i=1,...,n$. The differential $m_n$ of
the chain complex $M(n)$ is much simpler. Namely, the only non-trivial differential is
$m_n(\zzz_{i+1})=\yyy_{i}$ for $i=1,...,n-1$.\\

The map
$\Psi_1(n)$  will replace the intersection point $s$ with
the intersection point $x$, and will replace the intersection point $t$ with the intersection
point $z_1$. In other words, $\Psi_1(n)(\sss_i)=\xxx_i$ for $i=1,...,n-1$
and $\Psi_1(n)(\ppp)=\zzz_1$, while $\Psi_1(n)$ is trivial on the other
generators. Next, note that the map $\Phi(n)$ changes
$x$ to $r$. So we should have $\Phi(n)(\xxx_i)=\rrr_i$ for $i=1,...,n-1$.
It is also easy to see that $\Psi_3(n)$ would replace $s$ with $y$, which
means that the map is given by $\Psi_3(n)(\sss_i)=\yyy_i$ for $i=1,...,n-1$.\\

There is one domain in this Heegaard diagram which is a  hexagon and
may be used together with the triangle corresponding to
the map $\Psi_2(n)$, and some of the rectangles, to build up a
pentagon with vertices $z,z_i,w_i,p$ and $s$ (note that
despite the existence of this hexagon, the differentials may yet
be described combinatorially). It is not hard to see  that
this pentagon contributes to the triangle map (see \cite{Sar}).
Correspondingly, we obtain that $\Psi_2(n)(\sss_i)=(\zzz_i)$,
$i=1,...,n-1$. There are also
contributions from small triangle with vertices $z,r$ and $y$, and also
from the triangle with vertices $z,t$ and $z_n$. These correspond to
$\Psi_2(n)(\rrr_i)=\yyy_i,\ i=1,...,n-1$ and $\Psi_2(n)(\ppp)=\zzz_n$.
Putting all these together finishes our study of the chain complexes
and the maps associated with the Heegaard diagram $H_n$.\\

If $M=M(K)$ and $L=L(K)$ are the chain complexes associated
with the knot $K$ (and come from Heegaard diagrams of the
type discussed earlier), for a generator $\xxx$ of the chain complex
$M(n)$, let $M[\xxx]$ denotes the the copy of $M$ associated with  $\xxx$
(and define $L[\yyy]$ for a generator $\yyy$ of $L(n)$  similarly).
Some simple algebra implies that if we remove both
$L[\zzz_{i+1}]$ and $L[\yyy_{i}]$ from the complex, the homology of
the new complex is the same as the homology of the initial
complex). Once we do these simplifications, the remaining parts of the
complex will take the following form:
\begin{displaymath}
\begin{diagram}
L[\ppp]&\rTo^{i_L}&L[\rrr_{n-1}]&&&&L[\rrr_{n-2}]&&...&...&&L[\rrr_1]\\
       &&\uTo^{\Phi}&&&\ruTo&\uTo^{\Phi}&&&&\ruTo^{i_L}&\uTo^{\Phi}\\
\dTo^{\Psi_2}&&M[\xxx_{n-1}]&&\ldLine^{i_L}&&M[\xxx_{n-2}]&&...&...&&M[\xxx_1]\\
&&\uTo{\Psi_2} &\ruLine&&&\uTo{\Psi_2} &\ruTo^{i_L}&&&&    \\
&&L[\sss_{n-1}]&&&&L[\sss_{n-2}]&&...&...&&\uTo^{\Psi_2}\\
M[\zzz_1]&&&&&\lTo_{\Psi_1}&&&&&&L[\sss_1].
\end{diagram}
\end{displaymath}

\section{Algebraic simplifications}
Let $C_\bullet(K)$ for $\bullet\in\{\infty,1,0\}$ and the maps
$\phi,\phibar,\psi$ and $\psibar$ be defined as in \cite{Ef-splicing}
(or as described earlier in the introduction).
It was shown in \cite{Ef-splicing} that the complex $L(K)$ may be realized as the mapping
cone of the chain map $\phi$ and that $M(K)$ may be realized as the
mapping cone of the chain map $\psi$, for an appropriate choice of
the Heegaard diagram. Furthermore, the map $\Phi:M(K)\ra L(K)$
will be described as the identity map of $C_1$. The map $\Psi_1:L(K)\ra M(K)$
will be the map $\phibar:C_\infty \ra C_1$, and the map
$\Psi_2:L(K)\ra M(K)$ is the map $\psibar:C_1\ra C_0$. Using this description,
we may understand the blocks in the above diagram which are of the form
$$\xra{\Psi_2}M[\xxx_i] \xra{\Phi} L[\rrr_i]\xleftarrow{i_L} L[\sss_{i+1}]
\xra{\Psi_2}.$$
Doing the above substitutions will replace the above block with the
following
\begin{displaymath}
\begin{diagram}
...&&C_1(K)&&C_\infty(K)&\lTo^{i_\infty}&C_\infty(K)&&...\\
&&\dTo{\psi}&\rdTo{i_1}&\dTo{\phi}&&\dTo{\phi}&&\\
...&\rTo_{\psibar}&C_0(K)&&C_1(K)&\lTo_{i_1}&C_1(K)&\rTo_{\psibar}&...,
\end{diagram}
\end{displaymath}
where $i_\bullet$ is the identity map of $C_\bullet(K)$. When we compute the
homology of the big complex, the two $C_\infty(K)$ terms in this block may
be canceled against each-other, since no differential goes to
the block $C_\infty(K)\xleftarrow{i_\infty} C_\infty(K)$. In the remaining
part of this block we may replace
$$...\xleftarrow{\psi}C_1(K)\xra{i_1}C_1(K)\xleftarrow{i_1}C_1(K)
\xra{\psibar} ...$$
by $...\xleftarrow{\psi}C_1(K)\xra{\psibar}...$ without changing the homology.
There are $n-1$ such blocks in the above diagram, which are followed by
one-another.
The remaining part of the above diagram is a sequence of the form
$$
\xra{\Psi_2}M[\zzz_1] \xla{\Psi_1}L[\sss_1]\xra{\Psi_2} .$$
In terms of the complexes $C_\bullet(K)$, this part of the
diagram takes the form
\begin{displaymath}
\begin{diagram}
...&\rTo{\psibar}&C_0(K)&\lTo{\psi}&C_1(K)&\lTo{\phibar}&C_\infty(K)&\rTo{\phi}&C_1(K)&\rTo{\psibar}&....
\end{diagram}
\end{displaymath}
Putting all these pieces together, we obtain a much simpler diagram.
As in \cite{Ef-splicing}, the complexes
$C_\bullet(K)$ may be replaced by their homology groups $\Hbb_\bullet=
\Hbb_\bullet(K)$, since $\psi\circ \phibar=\psibar\circ \phi=0$ and the induced
maps
\begin{displaymath}
\begin{split}
&\theta:\text{Ker}(\phi_*)\subset \Hbb_\infty\ra \text{Coker}(\psibar_*)=\Hbb_0/\text{Im}(\psibar_*),\ \text{and}\\
&\thetabar:\text{Ker}(\phibar_*)\subset \Hbb_\infty\ra \text{Coker}(\psi_*)=\Hbb_0/\text{Im}(\psi_*)
\end{split}
\end{displaymath}
are both trivial. We may collect all the above observations to
 give the following surgery formula:
\begin{thm}
Let $K$ be a knot in a homology sphere $Y$, and let the complexes $\Hbb_\bullet$, $\bullet
\in \{\infty,1,0\}$ and the maps $\phi,\phibar,\psi,\psibar$ between them be as before.
The homology of $Y_n(K)$, the manifold obtained by $n$-surgery on $K$, may be obtained
as the homology of the following complex
\begin{displaymath}
\begin{diagram}
\Hbb_0&\lTo^{\psi}&\Hbb_1&\rTo^{\psibar}&\Hbb_0&\lTo^{\psi}&...&\rTo^{\psibar}&\Hbb_0\\
     &\luTo_{\psibar}&   &              &      &           &   & \ruTo_{\psi} &    \\
&  &\Hbb_1&\rTo_{\phi}&\Hbb_\infty&\lTo_{\phibar}&\Hbb_1&    &    \\
\end{diagram}
\end{displaymath}
where the total number of terms in the first row is $2n-1$.
\end{thm}

With a completely similar technique, and a little more computations, one can show the
rational surgery formula that was stated in theorem~\ref{thm:main}.


\end{document}